\documentclass{article}
\begin{document}
\centerline{\large\bf On the numbers that are sums of three  cubes}

\[
\]

\centerline{N.D. Bagis}
\centerline{Aristotele University of Thessaloniki-AUTH}
\centerline{Thessaloniki-Greece}
\centerline{e-mail:nikosbagis@hotmail.gr}

\begin{quote}

\[
\]

\centerline{\bf Abstract}
We examine what integers are representable as sums of three cubes. We also provide formulas for the number of representations of $x^3+y^3+z^3=n$ under the condition  $x+y+z=t$. Also we show how the problem of three cubes is related to $abc-$conjecture. 
      
\[
\]

\textbf{keywords}: \textrm{Sums of three cubes; Diophantine equations; Higher forms; Representations of integers}

\end{quote}

\section{Introduction.}

We will proceed to find all integer solutions $(x,y,z)$ of 
\begin{equation}
x^3+y^3+z^3=n
\end{equation}
with height
\begin{equation}
x+y+z=t.
\end{equation}
The above problem (without condition (2)) is rather difficult and is open for more than hundred years years. In recent days, with the help of computers many unknown results come to light. For example it simple to state that 
$$
1^3+1^3+1^3=4^3+4^3+(-5)^3=3.
$$
But it is very difficult to find that
$$
569936821221962380720^3+(-569936821113563493509)^3+
$$
$$
+(-472715493453327032)^3=3.
$$
In 1992, Roger Heath-Brown conjectured that every $n$ unequal to 4 or 5 modulo 9 has infinitely many representations as sums of three cubes. This conjecture until now days remains unproved.

\section{Analyzing the problem.}

Set
\begin{equation}
x^2+y^2+z^2=n_1
\end{equation}
The cases (1) together with (2) are solvable with Wolfram alpha.\\

Set
\begin{equation}
x+y+z=t\textrm{, }x y+y z+z x=r\textrm{, }x y z=s.
\end{equation}
Then
\begin{equation}
x^2+y^2+z^2=t^2-2r=n_1.
\end{equation}
Hence 
\begin{equation}
t=\pm\sqrt{n_1+2r}
\end{equation}
Also then
\begin{equation}
x^3+y^3+z^3=n\Leftrightarrow 3s-3rt+t^3=n.
\end{equation}
Setting the value of (6) to (7) we get
\begin{equation}
n\pm(r-n_1)\sqrt{2r+n_1}=3s.
\end{equation}
We set also $p$ be such that
\begin{equation}
p=3s-n.
\end{equation}
Equation (7) becomes finally
\begin{equation}
t(t^2-3r)=-p,
\end{equation}
which is solvable (however this means nothing since for solving (10) we need to know where $p$ exactly varies. Clearly we have
\begin{equation}
|xy+yz+zx|\leq x^2+y^2+z^2\Leftrightarrow |r|\leq n_1
\end{equation}
Also hold the following identity of Euler
\begin{equation}
x^3+y^3+z^3-3xyz=(x+y+z)(x^2+y^2+z^2-x y-y z-z x)
\end{equation}
Hence
$$
n-3s=t(n_1-r)\Leftrightarrow-p=t(n_1-r)\Leftrightarrow \frac{p}{t}+n_1=r\Rightarrow\left|\frac{p}{t}+n_1\right|=|r|\leq n_1
$$
Hence
\begin{equation}
-2n_1\leq\frac{p}{t}\leq 0
\end{equation}
Also from (6) we have 
\begin{equation}
t^2=n_1+2r\Leftrightarrow \frac{t^2-n_1}{2}=r\Rightarrow \left|\frac{t^2-n_1}{2}\right|=|r|\leq n_1\Rightarrow-2n_1\leq t^2-n_1\leq 2n_1
\end{equation}
Hence $0\leq t^2\leq 3n_1$ and $r\geq -\frac{n_1}{2}\Rightarrow \frac{p}{t}\geq -n_1-\frac{n_1}{2}\Rightarrow0\geq\frac{p}{t}\geq-\frac{3n_1}{2}$. Also from (7)\\ 
$$
n-3s=t(n_1-r)
$$
Hence given 
$$
x^3+y^3+z^3=n,
$$
with height 
$$
x+y+z=t
$$
If we set 
$$
\epsilon=n_1-r=|n_1-r|\geq 0\textrm{ (always)},
$$
then from (7) we get
$$
n-3s=t\epsilon\Leftrightarrow \epsilon=t^2-3r
$$
and we have to solve the system
\begin{equation}
3tr-3s=t^3-n\textrm{, }n_1+2r=t^2.
\end{equation}
with respect to $r,s,n_1$. 
We have solution if $\textrm{gcd}(3t,3)=3|(t^3-n)$
\begin{equation}
r=C_1\textrm{, }s=tC_1+\frac{n-t^3}{3}\textrm{, }n_1=-2C_1+t^2\textrm{, }C_1\in\textbf{Z}.
\end{equation}
Now consider the product
$$
(x-t)(y-t)(z-t)=xyz-t(x y+y z+z x)+t^2(x+y+z)-t^3=
$$
$$
=x y z-t(x y+y z+z x)=s-t r=
$$
$$
=tC_1+\frac{n-t^3}{3}-t C_1=\frac{n-t^3}{3}
$$
Hence if we solve the equation
\begin{equation}
(x-t)(y-t)(z-t)=\frac{n-t^3}{3},
\end{equation}
we are done. Set $x=A+t$, $y=B+t$, $z=C+t$ and the equation (17) becomes
\begin{equation}
ABC=\frac{n-t^3}{3}.
\end{equation}
The general equation $ABC=N$ have $L(N)=4\sum_{0<d|N}\tau(d)$ solutions (here $\tau(n)=\sum_{d|n}1$ is the number of positive divisors of $n$). Hence when $N=N_1=\frac{n-t^3}{3}$, we have
\begin{equation}
L(N_1)=\sum_{ABC=(n-t^3/3)}1=2\sum^{*}_{\textrm{\scriptsize abs}(d)|N_1}\tau(d)=4\sum_{0<d|N_1}\tau(d),
\end{equation}
solutions (where the asterisk in the summation means that we exclude the case $d=0$).\\By (18),(17),(16),(15), we have\\
\\
\textbf{Theorem 1.}\\
The equation $x^3+y^3+z^3=n$, with height $x+y+z=t$, $n\neq t^3$, have finite number of solutions. Moreover the problem (1),(2) is equivalent to
\begin{equation}
ABC=\frac{n-t^3}{3}\textrm{, }A+B+C=-2t.
\end{equation}
Actually $A=x-t$, $B=y-t$, $C=z-t$, $t=x+y+z$.\\
\\

Hence given $n,t$, equation (18) have $L(N_1)$ solutions, where 
\begin{equation}
L(N):=4\sum_{d|N}\tau(d)
\end{equation}
and $N_1=\left|\frac{n-t^3}{3}\right|$. Taking one by one the $L(N_1)$ solutions of (18), we evaluate all cases such
$$
A+B+C=-2t
$$
and solve (1),(2) according  to Theorem 1.\\
\\
Also we have the next\\
\\
\textbf{Theorem 2.}\\ Given integers $n,t$, with $n\neq t^3$, the number of representations of the integer $n$ as   
\begin{equation}
x^3+y^3+z^3=n\textrm{, }(x,y,z)\in\textbf{Z}^3,
\end{equation}
with height
\begin{equation}
x+y+z=t,
\end{equation}
is
\begin{equation}
R(t,n)=\sum_{\scriptsize
\begin{array}{cc}
0\neq\textrm{\scriptsize abs\normalsize}(d)|\textrm{\scriptsize abs\normalsize}(N(t,n))\\
0\neq\textrm{\scriptsize abs\normalsize}(\delta)|\textrm{\scriptsize abs\normalsize}(d)\\
\delta+d/\delta+N(t,n)/d=-2t	
\end{array}\normalsize}1,
\end{equation}
where $N(t,n)=\frac{n-t^3}{3}$. Note that if $n-t^3\neq 0(\textrm{mod} 3)$, we don't have representations. Also if $n=t^3$, we have infinite representations. In this case (and for convenience) we set in 0 in (24).\\
\\
\textbf{Theorem 3.}\\
Assume that $n$ is a fixed integer, then the number of representations of $n$ in the form 
\begin{equation}
x^3+y^3+z^3=n\textrm{, }(x,y,z)\in\textbf{Z}^3,
\end{equation}
with $-j\leq x+y+z\leq j$, where $j\geq 0$, is
\begin{equation}
C(j,n)=\sum^{j}_{t=-j}R(t,n),
\end{equation}
where $R(t,n)$ is that of (24). Moreover if $n_0$ is an integer such that $n_0-t^3\neq0$ and $n_0-t^3\equiv0(\textrm{mod}3)$. Then the total number of representations of $n_0$ in (25) is $\lim_{j\rightarrow+\infty}C(j,n_0).$\\
\\
\textbf{Note.} Assume that exists $n_0,t_0\in\textbf{N}$ such that $R(t,n_0)=0$, for all $t$ with $|t|\geq t_0>0$ and for at least one $t_1<t_0$, we have $R(t_1,n_0)\geq 1$.\\ 
\textbf{i)} If $n_0\equiv0(\textrm{mod}3)$ and $t_1\equiv 1,2(\textrm{mod}3)$, then $R(t_1,n_0)=0$, which is not true. Hence $t_1\equiv 0(\textrm{mod}3)$.\\
\textbf{ii)} If $n_0\equiv 1(\textrm{mod}3)$, then if $t_1\equiv0,2(\textrm{mod}3)$, we get $R(t_1,n_0)=0$, which is not true. Hence $t_1\equiv 1(\textrm{mod}3)$.\\In the same way\\
\textbf{iii)} If $n_0\equiv 2(\textrm{mod}3)$, then $t_1\equiv 2(\textrm{mod}3)$.\\Hence in general we can write
\begin{equation}
t_1\equiv n_0(\textrm{mod}3).
\end{equation}
\\

If $N_1=\frac{n-t^3}{3}$ and $n-t^3\equiv0(\textrm{mod}3)$, then assuming that $d$ runs through all integer divisors of $N_1$ and $\delta$ is divisor of $d$ and hence of $N_1$, in order to have solution, we must have    
\begin{equation}
\delta+d/\delta+N_1/d=-2t.
\end{equation}
Solving with respect to $\delta$ the above equation, we get 
\begin{equation}
\delta=-\left[\frac{N_1}{2d}+t\pm\sqrt{-d+\left(\frac{N_1}{2d}+t\right)^2}\right].
\end{equation}
Hence $N_1$ must be even. However exist cases such $N_1$ is even and we have no solution i.e. $t=-2$ and $n=4$, $N_1=4$. Hence when $N_1$ is odd we have no solutions. In general relation (29) give us the following\\
\\
\textbf{Theorem 4.}\\
Given $n$, $t$ integers, we set $N_1:=\frac{n-t^3}{3}$. Then the number of representations in the form 
\begin{equation}
x^3+y^3+z^3=n,
\end{equation}
with height
\begin{equation}
x+y+z=t,
\end{equation}
is
\begin{equation}
R(t,n)=-\sum_{\scriptsize
\begin{array}{cc}
	0<d|N_1\\
	(N_1/d+2t)^2=4d
\end{array}
\normalsize}1+2\sum_{
	0\neq\textrm{\scriptsize abs\normalsize}(d)|N_1}S_{\textbf{\scriptsize Z\normalsize}_{+}}\left(\left(\frac{N_1}{d}+2t\right)^2-4d\right),
\end{equation}
where $S_{\textbf{\scriptsize Z\normalsize}_{+}}(n)$ is 1 if $n$ is a nonnegative square of integer and 0 otherwise.\\ 
\\
\textbf{Proof.}\\
In order to have a representation added in the sum (24), a divisor $d$ of $N_1=\frac{n-t^3}{3}$ must be such  
\begin{equation}
x_d\pm\sqrt{x_d^2-d}\in\textbf{Z},
\end{equation}
where $x_d=\frac{N_1}{2d}+t$. Hence the number of representations becomes
\begin{equation}
R(t,n)=\sum_{\scriptsize
\begin{array}{cc}
	0\neq\textrm{\scriptsize abs\normalsize}(d)|N_1\\
	x_d\pm\sqrt{x_d^2-d}\in\textbf{\scriptsize Z\normalsize}
\end{array}
\normalsize}1-\sum_{\scriptsize
\begin{array}{cc}
	0\neq\textrm{\scriptsize abs\normalsize}(d)|N_1\\
  x_d^2=d
\end{array}
\normalsize}1.
\end{equation}
Here we take the negative part to avoid zero discriminants. Hence if $u=x_d-\sqrt{x_d^2-d}$ and $v=x_d+\sqrt{x_d^2-d}$, then $uv=d$ and $u+v=2x_d$ and we also have
\begin{equation}
R(t,n)=3\sum_{\scriptsize
\begin{array}{cc}
	0\neq\textrm{\scriptsize abs\normalsize}(d)|N_1\\
	0<|u|\leq d\\
	0<|v|\leq d\\
	uv=d\\
	u+v=N_1/d+2t
\end{array}
\normalsize}1.
\end{equation}
However the equivalence (28)$\Leftrightarrow$(29) is that which reduces the problem.\\  
\\
\textbf{Lemma 1.}\\
\textbf{1)} If $n=9k$, then\\ 
\textbf{i)} If $k=4k_1$ in order to have solutions of (1), we must have $t=6l$.\\ 
\textbf{ii)} If $k=4k_1+1$ then $t=6l+3$.\\
\textbf{iii)} If $k=4k_1+2$, then $t=6l$.\\ 
\textbf{iv)} If $k=4k_1+3$, then $t=6l+3$.\\
\textbf{2)} If $n=9k+1$, then\\
\textbf{i)} If $k=4k_1$, then $t=6l+1$.\\
\textbf{ii)} If $k=4k_1+1$, then $t=6l+4$.\\ 
\textbf{iii)} If $k=4k_1+2$, then $t=6l+1$.\\
\textbf{iv)} If $k=4k_1+3$, then $t=6l+4$.\\
\textbf{3)} If $n=9k+2$, then\\
\textbf{i)} If $k=4k_1$, then $t=6l+2$.\\
\textbf{ii)} If $k=4k_1+1$, then $t=6l+5$.\\ 
\textbf{iii)} If $k=4k_1+2$, then $t=6l+2$.\\
\textbf{iv)} If $k=4k_1+3$, then $t=6l+5$.\\
\textbf{4)} If $n=9k+3$, then\\
\textbf{i)} If $k=4k_1$, then $t=6l+3$.\\
\textbf{ii)} If $k=4k_1+1$, then $t=6l$.\\
\textbf{iii)} If $k=4k_1+2$, then $t=6l+3$.\\
\textbf{iv)} If $k=4k_1+3$, then $t=6l$.\\ 
\textbf{5)} If $n=9k+6$, then\\
\textbf{i)} If $k=4k_1$, then $t=6l$.\\
\textbf{ii)} If $k=4k_1+1$, then $t=6l+3$.\\
\textbf{iii)} If $k=4k_1+2$, then $t=6l$.\\
\textbf{iv)} If $k=4k_1+3$, then $t=6l+3$.\\ 
\textbf{6)} If $n=9k+7$, then\\
\textbf{i)} If $k=4k_1$, then $t=6l+1$.\\
\textbf{ii)} If $k=4k_1+1$, then $t=6l+4$.\\
\textbf{iii)} If $k=4k_1+2$, then $t=6l+1$.\\
\textbf{iv)} If $k=4k_1+3$, then $t=6l+4$.\\ 
\textbf{7)} If $n=9k+8$, then\\
\textbf{i)} If $k=4k_1$, then $t=6l+2$.\\
\textbf{ii)} If $k=4k_1+1$, then $t=6l+5$.\\
\textbf{iii)} If $k=4k_1+2$, then $t=6l+2$.\\
\textbf{iv)} If $k=4k_1+3$, then $t=6l+5$.\\ 
\\
\textbf{Remarks.} Actually it is known that when $n=9k\pm4$, then equation (1) have no representations at all.\\
\\
\textbf{Proof of Lemma 1.}\\
\textbf{1)} For $n=9k$ we have\\
\textbf{i)} Assuming $n=9k$, $k=4k_1$, then using $3|(n-t^3)$ and $N_1=\frac{n-t^3}{3}=even$ we get easily the first result.\\ 
\textbf{ii)} For $n=9k$, $k=4k_1+1$, then using $3|(n-t^3)$ and $N_1=\frac{n-t^3}{3}=even$ we get the second result.\\
\textbf{iii)} The third result ($n=9k$, $k=4k_1+2$) follows in the same way $\ldots$etc.\\
\\
\textbf{Theorem 5.}\\
When $n\equiv\pm4(\textrm{mod}9)$ the equation (1) have no solutions.\\
\\
\textbf{Proof.}\\
From Theorem 4 relation (32) we have (in order to have solutions)
$$
\left(\frac{N_1}{d}+2t\right)^2-4d=k^2,\eqno{(a)}
$$  
where $d$ is any divisor of $N_1$. From relation ($a$) we have that exists integer $y$ such that
$$
4d=y^2-k^2.\eqno{(b)}
$$ 
Clearly $k,y$ are both even or both odd.\\
\textbf{1.}\\
\textbf{i)} Assume that both $k,y$ are odd, then $k=2v+1$, $y=2u+1$. Then $4d=(4u^2+4u+1)-(4v^2+4v+1)=4(u^2+u-(v^2+v))$. Hence 
$$
d=u^2+u-v^2-v.
$$
Relation ($a$) can be written in the form
$$
N_1=yd-2td\Leftrightarrow n-t^3=3yd-6td\Leftrightarrow
$$
$$ 
n=3(2u+1)(u^2+u-v^2-v)-6t(u^2+u-v^2-v)+t^3.\eqno{(c)}
$$ 
From Lemma 1 we have if $n=9s+4$, $s$ even, then $t=4+6l$. Hence $n=4+18s_1\equiv4(\textrm{mod}18)$. Also then $t^3\equiv10(\textrm{mod}18)$, $6t\equiv6(\textrm{mod}18)$. Then relation ($c$) becomes
$$
n\equiv 3(2u+1)(u^2+u-v^2-v)-6(u^2+u-v^2-v)+10(\textrm{mod}18).
$$
As one can see by taking all the possible cases $u\equiv0,1,\ldots,17(\textrm{mod}18)$ and $v\equiv0,1,\ldots,17(\textrm{mod}18)$, that then $n\equiv 10,16(\textrm{mod}18)$. Hence $n$ is not $4(\textrm{mod}18)$.\\
In case $y=2u$, $k=2v$, then $d=u^2-v^2$ and
$$
n=3(u^2-v^2)(2u)-6t(u^2-v^2)+t^3.\eqno{(d)}
$$ 
This case lead us again to $n\equiv10,16(\textrm{mod}18)$, which is not of the form $n\equiv 4(\textrm{mod}18)$.\\
\textbf{ii)} From Lemma 1, the case $n=9s+4$, $s$ odd is $n\equiv13(\textrm{mod}18)$ and $t\equiv 1(\textrm{mod}6)$. Hence $t^3\equiv 1(\textrm{mod}18)$ and $6t\equiv 6(\textrm{mod}18)$. Hence then (when $y=2u+1$, $k=2v+1$)
$$
n\equiv 3(u^2+u-v^2-v)(2u+1)-6(u^2+u-v^2-v)+1(\textrm{mod}18).
$$ 
Taking all the possible cases we see that $n\equiv1,7(\textrm{mod}18)$, which is not of the form $n\equiv 13(\textrm{mod}18)$.\\
For $y=2u$, $k=2v$, we have
$$
n=3(u^2-v^2)(2u)-6(u^2-v^2)+1(\textrm{mod}18)\equiv1,7(\textrm{mod}18).
$$
Hence again is not proper.\\
\textbf{2.}\\
\textbf{i)} For the case $n=-4+9s$, $s$ even, we have $t=5+6l$. Hence $n\equiv-4(\textrm{mod}18)\equiv14(\textrm{mod}18)$, we have $t\equiv5(\textrm{mod}6)$. Hence $6t\equiv12(\textrm{mod}18)$ and $t^3\equiv17(\textrm{mod}18)$.\\Hence when $y=2u+1$, $k=2v+1$, then
$$
n\equiv3(u^2+u-v^2-v)(2u+1)-12(u^2+u-v^2-v)+17(\textrm{mod}18)
$$ 
and this leads to $n\equiv 11,17(\textrm{mod}18)$, which not proper.\\
For the case $y=2u$, $k=2v$, we have 
$$
n\equiv3(u^2-v^2)(2u)-12(u^2-v^2)+17(\textrm{mod}18)\equiv11,17(\textrm{mod}18)
$$
which is not proper.\\
\textbf{ii)} The case $n=-4+9s$, with $s$ odd give us $n\equiv5(\textrm{mod}18)$, $t\equiv2(\textrm{mod}6)$, $6t\equiv12(\textrm{mod}18)$, $t^3\equiv8(\textrm{mod}18)$, lead us again to no solution.\\
\\
Hence in any case the theorem is proved.\\
\\ 

From the above theorem we get the equivalent condition for to have solutions of (1),(2):\\
\\
\textbf{Theorem 6.}\\
If $N_1=\frac{n-t^3}{3}$ and exists $u,v$ integers such that 
\begin{equation}
N_1=(2u+1)(u^2+u-v^2-v)-2t(u^2+u-v^2-v)
\end{equation}
or
\begin{equation}
N_1=2u(u^2-v^2)-2t(u^2-v^2)
\end{equation}
then we have a solution of (1),(2) ($+1$ or $+2$ added in  $R(t,n)$). In case that both (36),(37) are not satisfied then $R(t,n)=0$.\\ 
i) In case of (37), with $k=0$ then we have at least one solution ($+1$ in $R(t,n)$).\\ 
ii) If $k\neq0$, then for each $(u,v)$ satisfying (37) we add $+2$ in $R(t,n)$.\\
iii) For each $(u,v)$ satisfying (36) we have $+2$ in $R(t,n)$.\\
iv) In any case in order to have $+1$ or $+2$ in $R(t,n)$, $N_1$ must be even and if $d|N_1$, then $d$ must be of the form $4d=y^2-k^2$, with $|y|=\left|N_1/d+2t\right|$.\\
\\
\textbf{Corollary 1.}\\
$R(0,n)>0$, whenever $n\equiv 0(\textrm{mod}6)$ and exists divisor $d$ of $n/3$ such that 
\begin{equation}
\frac{n^2}{9d^2}-4d=\textrm{perfect square}.
\end{equation}
\\

Assume now the equation
\begin{equation}  
9d^2(k^2+4d)=n^2
\end{equation}
where $n,k,d$ are in general integers. Assume $L(n)$ denotes the number of solutions of (39) for a given $n$. If we where able to prove that exists sequence of natural numbers $n_m$, $m=1,2,\ldots$, such that $n_m$ were increasing with $\lim_{m\rightarrow+\infty}L(n_m)=+\infty$ i.e. the number of solutions of (39) for $n=n_m$ will be growing without a bound as $m$ goes to infinity. Then in view of Theorems 2,3,4, the infinite numbers $n=n_m$ will correspond to
\begin{equation}
x^3+y^3+z^3=n_m\textrm{, }x+y+z=0\textrm{, }m=1,2,\ldots
\end{equation}
and have $L(n_m)=R(0,n_m)$ number of solutions and 
\begin{equation}
\lim_{m\rightarrow\infty} R(0,n_m)=+\infty.
\end{equation}
Hence the equations
\begin{equation}
x^3+y^3+z^3=n_m\textrm{, }m=1,2,\ldots,
\end{equation}
with $L^*(n_m)$ being the number of solutions of the $m-$th equation, will have at least finite number of solutions (from Theorem 3). But from (41) and Theorem 3 we will have  
\begin{equation}
\lim_{m\rightarrow+\infty}L^*(n_m)=+\infty.
\end{equation} 
Hence, as $m\rightarrow+\infty$, the number of solutions of (42) will then goes to infinity.\\
Numerically we can construct such numbers $n_m$. We can see also (numerically) that they satisfy the approximate relation
\begin{equation}
\sigma_1(n_m)\approx e^{\gamma}n_m\log\log n_m.
\end{equation}
Hence if we define the set
\begin{equation}
S(a,b)=\{n\in\textbf{N}:5040< a\leq n\leq b\textrm{ and }0.85<\frac{\sigma_1(n)}{e^{\gamma}n\log\log n}<1\}.
\end{equation}
Then $n=n_m\in(a,b)$ iff $n\in S(a,b)$, $L(n_{m+1})-L(n_m)=6$ and $n$ is the smallest possible i.e. from all $n\in\textbf{N}$ such that $L(n)=c=const>0$ ''say'' $c=24$, we choose the smallest one $n=19440$. We have calculated the first few values of $n_m$. These are $n_1=90$, $L(n_1)=12$; $n_2=720$, $L(n_2)=18$ $n_3=19440$, $L(n_3)=24$; $n_4=55440$, $L(n_4)=30$; $n_5=443520$, $L(n_5)=36$;$\ldots$. Note that Riemann's hypothesis is equivalent with the Robin$-$Ramanujan inequality
\begin{equation}
\sigma_1(n)<e^{\gamma}n\log\log n\textrm{, }\forall n>5040,
\end{equation}
where $\sigma_{\nu}(n)=\sum_{d|n}d^{\nu}$ is the sum of the $\nu-$th power of positive divisors of $n$. But of course this is no proof that such $n_m$ are infinite. In the next section we provide a proof based in $abc-$conjecture see Theorem 9 below.\\

For to solve (39), we have $n=3\cdot n_1$, where $n_1=2^{a_1}3^{a_2}\ldots p_s^{a_s}$ is the prime decomposition of $n_1$. Then $d=2^{b_1}3^{b_2}\ldots p_s^{b_s}$, with $c_i=a_i-b_i\geq 0$, $i=1,2,\ldots,s$ is any divisor of $n_1$ and (39) becomes
\begin{equation}
l^2-k^2=4d\textrm{, }\forall d|n_1.
\end{equation} 
This equation have classically $2\tau(d)$ number of solutions. But for these solutions we must have the extra condition: $l=2^{c_1}3^{c_2}\ldots p_s^{c_s}$ or equivalently $l=n/(3d)$.\\

In case $t=t_0\neq 0$, the equation 
$$
x^3+y^3+z^3=n\textrm{, }x+y+z=t_0,
$$
is equivalent to
$$
\left(\frac{n-t_0^3}{3d}+2t_0\right)^2-4d=k^2,
$$
where $d$ is divisor of $N_1=\frac{n-t_0^3}{3}$. This becomes 
\begin{equation}
l^2-k^2=4d\textrm{, }(l-2t_0)d=N_1.
\end{equation}   
 
\section{The reduced problem}

Assume the equation
\begin{equation}
xy(x-y-1)=n,
\end{equation}
where $x,y,n\in\textbf{Z}$. I rewrite equation (49) in the form
$$
(-x)y(-x-y-1)=n\Leftrightarrow xy(-x-y-1)=-n.
$$
Hence setting $z=-x-y-1$, we get that (49) is equivalent to
\begin{equation}
xyz=-n\textrm{, }x+y+z=-1.
\end{equation}
Now given $n$ integer, according to Theorem 2, equation (49) and hence (50) have
\begin{equation}
L(n)=R\left(\frac{1}{2},\frac{1-24n}{8}\right)=\sum_{\scriptsize
\begin{array}{cc}
0\neq\textrm{\scriptsize abs\normalsize}(d)|n\\
0\neq\textrm{\scriptsize abs\normalsize}(\delta)|\textrm{\scriptsize abs\normalsize}(d)\\
\delta+d/\delta-n/d=-1	
\end{array}\normalsize}1,
\end{equation}
number of solutions. I will show that for certain choices of $n$ we have always solution of (37). Set $n=pp_1(1+p+p_1)$, then exists always a divisor $d_0=pp_1>1$ of $n$ and a divisor $\delta_0=p$ of $d_0$, such that $\delta_0+d_0/\delta_0-n/d_0=-1$. Taking $n=n_{\nu}=p_{\nu}p_{\nu+1}(1+p_{\nu}+p_{\nu+1})$, where $p_{\nu}$ is the $\nu-$th prime, we have an infinite choice of positive integers $n=n_{\nu}$ such that $L(n_{\nu})>0$. Hence equation (49)$\Leftrightarrow$(50) have infinite number of $n$ with at least one solution.\\
\\
Also in view of Theorem 4 and relation (34) we get the next\\
\\
\textbf{Theorem 6.}\\
The equation (49)$\Leftrightarrow$(50) have number of solutions:
\begin{equation}
L(n)=-\sum_{\scriptsize
\begin{array}{cc}
	0<d|n\\
	(n/d-1)^2=4d
\end{array}\normalsize}1+2\sum_{0\neq \textrm{\scriptsize abs\normalsize}(d)|n
}S_{\textbf{\scriptsize Z\normalsize}_{+}}\left(\left(\frac{n}{d}-1\right)^2-4d\right)
\end{equation}
\\
\textbf{Proof.}\\
The equation (49)$\Leftrightarrow$(50) have
\begin{equation}
L(n)=\sum_{\scriptsize
\begin{array}{cc}
	0\neq \textrm{\scriptsize abs\normalsize}(d)|n\\
	y_d\pm\sqrt{y_d^2-d}\in\textbf{\scriptsize Z\normalsize}\\
	y_d^2>d
\end{array}
\normalsize}1-\sum_{\scriptsize
\begin{array}{cc}
	0\neq \textrm{\scriptsize abs\normalsize}(d)|n\\
	y_d^2=d,y_d\in\textbf{\scriptsize Z\normalsize}
\end{array}
\normalsize}1,
\end{equation}
where $y_d=\frac{n}{2d}-\frac{1}{2}$ number of solutions. Hence if $u=y_d-\sqrt{y_d^2-d}$ and $v=y_d+\sqrt{y_d^2-d}$, then $uv=d$ and $u+v=n/d-1$. Also then
\begin{equation}
L(n)=3\sum_{\scriptsize
\begin{array}{cc}
0\neq\textrm{\scriptsize abs\normalsize}(d)|n\\
0<|u|\leq d\\
0<|v|\leq d\\  
uv=d\\
u+v=n/d-1
\end{array}
\normalsize}1.
\end{equation}
\\

Assume now the form
\begin{equation}
x^2y+xy^2=n\textrm{, }(x,y)\in\textbf{Z}^2,
\end{equation}  
with $n\in\textbf{Z}$. This is equivalent to 
\begin{equation}
xyz=n\textrm{, }x+y-z=0\textrm{, }(x,y,z)\in\textbf{Z}^3
\end{equation}
and this to
\begin{equation}
xyz=-n\textrm{, }x+y+z=0\textrm{, }(x,y,z)\in\textbf{Z}^3.
\end{equation} 
If we manage to show that (57) have infinite solutions  $(x,y,z)\in\textbf{Z}^3$ for an infinite sequence $n=n_m\in\textbf{N}$, $m=1,2,\ldots$ and the number of solutions of each $n=n_m$ grows unboundedly. Then according to Theorem 3, equation $x^3+y^3+z^3=n$, will have also unboudeddly set of solutions. However (55)$\Leftrightarrow$(56)$\Leftrightarrow$(57). Assuming equation (55) with $x,y>0$, then according to $abc-$conjecture, we have that: for every $\epsilon>0$, exists a constant $K_{\epsilon}>0$ depending to $\epsilon$ only such that
\begin{equation}
\left(\textrm{rad}(n)\right)^{1+\epsilon}>\frac{1}{K_{\epsilon}}(x+y),
\end{equation}   
where
\begin{equation}
\overline{n}:=\textrm{rad}(n)=\prod_{\scriptsize
\begin{array}{cc}
	p|n\\
	p-prime
\end{array}\normalsize}p.
\end{equation}
Assuming  $Sol^{(+)}(n)=\{(x_1,y_1,z_1),(x_2,y_2,z_2),\ldots,(x_{\nu},y_{\nu},z_{\nu})\}$ is the set of positive solutions of (57) and $L^{(+)}(n)=\nu>0$ is the number of its solutions, then 
$$
L^{(+)}(n)\left(\textrm{rad}(n)\right)^{1+\epsilon}>\frac{1}{K_{\epsilon}}\sum^{\nu}_{k=1}z_k.
$$
Hence we have the next\\
\\
\textbf{Theorem 7.}\\
Equation (57)$\Leftrightarrow$(56)$\Leftrightarrow$(55) have for each $n$ positive integer 
\begin{equation}
\nu=L_0(n)=R(0,3n)=-\sum_{\scriptsize
\begin{array}{cc}
	0<d|n\\
	n^2=4d^3
\end{array}
\normalsize}1+2\sum_{0\neq\textrm{\scriptsize abs\normalsize}(d)|n}S_{\textbf{\scriptsize Z\normalsize}_{+}}\left(\frac{n^2}{d^2}-4d\right)
\end{equation}
solutions. Moreover in view of notes above we have 
\begin{equation}
\frac{1}{\nu}\sum^{\nu}_{k=1}z_k<K_{\epsilon}\overline{n}^{1+\epsilon}\textrm{, }\forall \epsilon>0.
\end{equation}
\\
\textbf{Theorem 8.}\\
The equation 
\begin{equation}
x^3+y^3+z^3=-3n\textrm{, }x+y+z=0\textrm{, }n\geq0
\end{equation}
have $L_0(n)$ solutions $(x,y,z)\in\textbf{Z}^3$. Under the $abc-$conjecture, for a given integer $n_0>0$ all solutions with $n\geq n_0$ have
\begin{equation}
|x|+|y|+|z|<C_{\epsilon}\overline{n}^{1+\epsilon},
\end{equation}
where $\epsilon>0$ and $C_{\epsilon}$ positive constant depending on the choice of $n_0$ and $\epsilon$.\\ 
\\
\textbf{Proof.}\\
Immediate consequence of Theorem 7, relation (58), the symmetry of (57) and Theorem 1.\\
\\

Continuing our arguments, we use the bounds conjectured in [1]:\\
\\
\textbf{Conjecture 1.}(Generalized $abc-$conjecture)\\
If $x,y,z$ are positive integers and $k=\textrm{rad}(n)$, where $x+y=z$, $n=xyz$, then there is a constant $C_1$ such that
\begin{equation}
z<k\exp\left[4\sqrt{\frac{3\log k}{\log\log k}}\left(1+\frac{\log\log\log k}{2\log\log k}+\frac{C_1}{\log\log k}\right)\right]
\end{equation} 
holds. Also there is a constant $C_2$ such that
\begin{equation}
z>k\exp\left[4\sqrt{\frac{3\log k}{\log\log k}}\left(1+\frac{\log\log\log k}{2\log\log k}+\frac{C_2}{\log\log k}\right)\right]
\end{equation}
holds infinetly often.\\
\\

From the above notes together with Theorems 1,3 we can show that\\
\\
\textbf{Theorem 9.}\\
Assuming the generalized $abc$ conjecture, the equation 
\begin{equation}
x^3+y^3+z^3=n\textrm{, }(x,y,z)\in\textbf{Z}^3,
\end{equation}
with $n=3b_m$, have unbounded set of solutions as $m\rightarrow+\infty$. The numbers $b_m$ are defined as follows: If $n_1=xy(x+y)$, $z=x+y$, $k=\textrm{rad}(xyz)$ and the triples $(x,y,z)\in\textbf{Z}^{{*}{3}}_{+}$ satisfy (65), then $b_m=n_1=xyz$.\\
\\  
\textbf{Proof.}\\
From Conjecture 1, since (65) is satisfied infinetly often, there exists certain sequence $b_m\in\textbf{N}$ such that  $n_1=b_m=xyz$, $m=1,2,\ldots$, $z=x+y$ of integers with $\lim b_m=+\infty$. Also $\lim\textrm{rad}(b_m)=\infty$. That is, from Cauchy inequality $z=x+y\geq 2\sqrt{xy}\Rightarrow z^3\geq 4b_m$. Using (64) we get $\textrm{rad}(b_m)\rightarrow+\infty$. Also $\lim L_0(b_m)=\infty$ i.e. (55) and hence (56)$\Leftrightarrow$(57) (with $n=n_1=b_m$) have unbounded set of solutions. This can be proved as follows: 
Assume the equation (66), with $n=b_m\in\textbf{Z}^{*}_{+}$. We denote with 
$$
Sol^{(+)}_{m}=\{(x_1,y_1,z_1),(x_2,y_2,z_2),\ldots,(x_{\nu},y_{\nu},z_{\nu})\},
$$ 
the set of positive integer triples $(x,y,z)=(x,y,x+y)$ who satisfy (56) when $n=b_m$. Then assume that as $m\rightarrow+\infty$ we have $\nu=\nu_{m}<\infty$ (is bounded). But  $k=\textrm{rad}(xyz)=\textrm{rad}(b_m)\rightarrow+\infty$ as $m\rightarrow\infty$. Also according to Conjecture 1 we can write
$$
\frac{1}{\nu}\sum^{\nu}_{i=1}z_{i}> k\exp\left[4\sqrt{\frac{3\log k}{\log\log k}}\left(1+\frac{\log\log\log k}{2\log\log k}+\frac{C_2}{\log\log k}\right)\right]\rightarrow+\infty.
$$
If $\nu=\nu_m$ were bounded then 
$$
\frac{1}{\nu}\sum^{\nu}_{i=1}z_{i}<\infty,
$$
which is impossible. Hence $\lim_{m\rightarrow\infty}\nu_{m}=+\infty$.\\
Hence as consequence of Theorems 1,2,3, we have for $n=3b_m$, that, the number of solutions of (66) is 
$$
\lim_{j\rightarrow+\infty}C(j,n)=\lim_{j\rightarrow+\infty}\sum_{j\geq |t|>0}R(t,3b_m)+R(0,3b_m)\geq R(0,3b_m).
$$ 
But $R(0,3b_m)$ is exactly the number of solutions of (57). Hence
$$
\lim_{m\rightarrow+\infty}R(0,3b_m)=+\infty.
$$
Hence when $m\rightarrow\infty$ the number of solutions of equation (66) with $n=n_m$ goes to infinity.\\
\\
\textbf{Theorem 10.}\\
Assume that $\nu=\nu_{m}$ is the number of solutions of 
\begin{equation}
xyz=b_m\textrm{, }z=x+y\textrm{, }x,y,z>0, 
\end{equation}
where $b_m$ is as defined in Theorem 9. Then under the generalized $abc-$conjecture we have $k=\textrm{rad}(xyz)\rightarrow+\infty$ as $m\rightarrow+\infty$ and
\begin{equation}
\lim_{m\rightarrow+\infty}\nu_{m}=+\infty.
\end{equation}
Further when $m\rightarrow+\infty$, we have the following asymptotic expansion
\begin{equation}
\frac{1}{\nu}\sum^{\nu}_{i=1}z_i=k\exp\left[4\sqrt{\frac{3\log k}{\log\log k}}\left(1+\frac{\log\log\log k}{2\log\log k}+O\left(\frac{1}{\log\log k}\right)\right)\right].
\end{equation}

\[
\]

\centerline{\textbf{References}}

[1]: Robert Olivier, Stewart Cameron L., Tenenbaum Gerald. ''A refinement of the abc conjecture''. Bulletin of the London Mathematical Society.\textbf{46}(6): 1156-1166.

\end{document}